\documentstyle[12pt,fleqn,amssymb,epsfig,a4]{article}

\newtheorem{theorem}{Theorem}[section]
\newtheorem{proposition}[theorem]{Proposition}
\newtheorem{lemma}[theorem]{Lemma}

\newtheorem{corollary}[theorem]{Corollary}
\newtheorem{claim}{Claim}[theorem]
\newtheorem{conjecture}[theorem]{Conjecture}

\begin{document}

\title{1-factor and cycle covers of cubic graphs}

\vspace{3cm}
  
\author{Eckhard Steffen\thanks{Paderborn Institute for Advanced Studies in 
		Computer Science and Engineering,
		Paderborn University,	
		Warburger Str. 100, 33098 Paderborn, Germany; 
		es@upb.de}}

\date{}

\maketitle

\begin{abstract}
{\small{Let $G$ be a bridgeless cubic graph. Consider a list of $k$ 1-factors of $G$. Let $E_i$ be the 
set of edges contained in precisely $i$ members of the $k$ 1-factors. Let $\mu_k(G)$ be the smallest $|E_0|$ over
all lists of $k$ 1-factors of $G$. 

Any list of three 1-factors induces a core of a cubic graph. We use results on the structure of cores to prove 
sufficient conditions for Berge-covers and for the existence of three 1-factors with empty intersection. 
Furthermore, if $\mu_3(G) \not = 0$, then $2 \mu_3(G)$ is an upper bound for the girth of $G$. 
We also prove some new upper bounds for the length of shortest cycle covers of bridgeless cubic graphs.

Cubic graphs with $\mu_4(G) = 0$ have a 4-cycle cover of length $\frac{4}{3} |E(G)|$
and a 5-cycle double cover. These graphs also satisfy two conjectures of Zhang \cite{Zhang_book}. 
We also give a negative answer to
a problem stated in \cite{Zhang_book}. 
}}
\end{abstract}  

\section[]{Introduction}

Graphs in this paper may contain multiple edges or loops. An edge is a loop if its two ends are incident to the same vertex. The degree 
of a vertex $v$ of a graph $G$ is the number of ends of edges that are incident to $v$.
A 1-factor of a graph $G$ is a spanning 1-regular subgraph of $G$.  Hence, a loop cannot be an edge of a 1-factor, and
a bridgeless cubic graph does not contain a loop. 
One of the first
theorems in graph theory, Petersen's Theorem from 1891 \cite{Petersen_1891}, states that every bridgeless cubic graph has a 1-factor.

Let $G$ be a cubic graph, $k \geq 1$, and $S_k$ be a list of $k$ 1-factors of $G$. By a list we mean a collection with possible repetition.
For $i \in \{0, \dots ,k\}$  let $E_i(S_k)$ 
be the set of edges which are in precisely $i$ elements of $S_k$. 
We define $\mu_k(G) = \min \{ |E_0(S_k)| \colon S_k \mbox{ is a list of $k$ 1-factors of } G\}$. If there is no harm of confusion we
write $E_i$ instead of $E_i(S_k)$.

If $\mu_3(G) = 0$, then $G$ is 3-edge-colorable. Cubic graphs with chromatic index 4 are the subject of many papers, 
since these graphs are potential counterexamples to many hard conjectures, see e.g.~\cite{Gunnar_2011_paper}.
Cyclically 4-edge-connected cubic graphs with chromatic index 4 and girth at least 5 are also called  {\em snarks}.

If $\mu_5(G)=0$, then the edges of $G$ can be covered by five 1-factors. This is conjectured to be true for all bridgeless cubic graphs
by Berge (unpublished, see e.g.~\cite{{Fouquet_09}, {Mazzuoccolo_2011}}). 

\begin{conjecture} [Berge Conjecture] \label{Berge_Conjecture}
Let $G$ be a cubic graph. If $G$ is bridgeless, then $\mu_5(G) = 0$.
\end{conjecture}

We say that $G$ has a Berge-cover, if it satisfies Conjecture \ref{Berge_Conjecture}.
Conjecture \ref{Berge_Conjecture} is true if the following conjecture is true, which is attributed to Berge in \cite{Seymour_79}. 
This conjecture was first published in a paper by Fulkerson \cite{Fulkerson_1971}.

\begin{conjecture} [Berge-Fulkerson Conjecture \cite{Fulkerson_1971}] \label{Fulkerson}
Let $G$ be a cubic graph. If $G$ is bridgeless, then $G$ has six 1-factors such that $E_2 = E(G)$. 
\end{conjecture}

We say that $G$ has a Fulkerson-cover, if it satisfies Conjecture \ref{Fulkerson}.
Mazzuoccolo  \cite{Mazzuoccolo_2011} proved that if  all bridgeless cubic graphs have a Berge-cover, 
then all bridgeless cubic graphs have a Fulkerson-cover. Hence, Conjectures \ref{Berge_Conjecture} and \ref{Fulkerson} are equivalent. 
Mazzuoccolo's result does not prove the equivalence of a Berge- and a Fulkerson-cover for a single bridgeless cubic graph. 
It is an open question whether a given cubic graph with a Berge-cover has a 
Fulkerson-cover as well.

The following conjecture of Fan and Raspaud is true if Conjecture \ref{Fulkerson} is true.
 
\begin{conjecture} [\cite{Fan_Raspaud_94}] \label{Fan_Raspaud}
Every bridgeless cubic graph has three 1-factors such that $E_3 = \emptyset$. 
\end{conjecture}

Eulerian graphs will be called {\em cycles} in the following. A {\em cycle cover}  of a graph $G$ is a set ${\cal C}$ of
cycles such that every edge of $G$ is contained in at least one cycle. It is a {\em double cycle cover}
if every edge is contained in precisely two cycles, and  a $k$-(double) cycle cover ($k \geq 1$) if ${\cal C}$
consists of at most $k$ cycles. Circuits of length 2 are allowed, but the two edges of such a circuit must be different. 
The following conjecture was stated by Celmins and Preissmann independently.

\begin{conjecture} [see \cite{Zhang_book}] \label{5-cdcc}
Every bridgeless graph has a 5-cycle double cover. 
\end{conjecture}

This conjecture is equivalent to its restriction to cubic graphs. 
In this case a cycle is a 2-regular graph. If it is connected, then we call it {\em circuit}; that is, a cycle is a set of disjoint circuits.
A cycle is {\em even} if it consists of even circuits. 
The length of a cycle cover ${\cal C}$ is the number of edges of ${\cal C}$. A 3-edge-colorable cubic graph $G$ has 
a 2-cycle cover of length $\frac{4}{3}|E(G)|$. Alon and Tarsi stated the following conjecture.

\begin{conjecture} [\cite{Alon_Tarsi}] \label{7/5 cover}
Every bridgeless graph $G$ has a cycle cover of length at most $\frac{7}{5}|E(G)|$.
\end{conjecture}

In Section \ref{gamma_3} we introduce the core of a cubic graph. 
Using structural properties of cores we show that if $\mu_3(G) \not = 0$, then
$2 \mu_3(G)$ is an upper bound for the girth of a cubic graph $G$. 
If $G$ is a bridgeless cubic graph without non-trivial 3-edge-cut and $\mu_3(G) \leq 4$, then $G$ has a Berge-cover.
If $G$ is a bridgeless cubic graph
and $\mu_3(G) \leq 6$, then $G$ satisfies Conjecture \ref{Fan_Raspaud}.  If $\mu_3(G) < girth(G)$, 
then $E_3= \emptyset$ for any three 1-factors with $|E_0| = \mu_3(G)$.
We  prove some new upper bounds for the length of a shortest cycle cover of bridgeless cubic graphs.
If $G$ is triangle-free and 
$\mu_3(G) \leq 5$, then $G$ has an even 3-cycle cover of length at most $\frac{4}{3}|E(G)| + 2$.
With those methods, some earlier results  of \cite{Fan_Raspaud_94} and \cite{Macajova_etal_2011}
are improved to even 3-cycle covers. 

Section \ref{gamma_4} proves that if $\mu_4(G) \in \{0,1,2,3\}$, then $G$ has 4-cycle cover of length 
$\frac{4}{3}|E(G)| + 4 \mu_4(G)$.
If $\mu_4(G)=0$, then there is a shortest 4-cycle cover of $G$ which is even and  and has edge depth at most 2,
i.e.~these graphs satisfy two conjectures of Zhang (Conjectures 8.11.5 and  8.11.6 in \cite{Zhang_book}).
Furthermore, $G$ has a 5-cycle double cover. We also give a negative answer to Problem 8.11.4 of \cite{Zhang_book}.

It seems that snarks which have only cycle covers with more than $\frac{4}{3}|E(G)|$ edges are rare.
There are only two such graphs with at most 36 vertices \cite{Gunnar_2011_paper}. One of these graphs is the 
Petersen graph and the other one has 34 vertices. Both graphs have a cycle cover of length $\frac{4}{3}|E(G)| + 1$. 

The paper concludes with a few remarks on hypohamiltonian snarks.

\section[]{Three 1-factors and the core of a cubic graph}  \label{gamma_3}

Let $M$ be a graph. If $X \subseteq E(M)$, then $M[X]$ denotes the graph whose vertex set consists of all vertices of edges
of $X$ and whose edge set is $X$. Let $A$, $B$ be two sets, then $A \oplus B$ denotes their symmetric difference.

Let $G$ be a cubic graph and $S_3$ be a list of three 1-factors $M_1, M_2, M_3$ of $G$. Let ${\cal M} = E_2 \cup E_3$, ${\cal U} = E_0$, and
$|{\cal U}| = k$.
The  {\em $k$-core} of $G$ with respect to $S_3$ (or to $M_1, M_2, M_3$)  is the subgraph $G_c$ of $G$ 
which is induced by ${\cal M} \cup {\cal U}$; that is, $G_c = G[{\cal M} \cup {\cal U}]$. 
If the value of $k$ is irrelevant, then we
say that $G_c$ is a core of $G$. 
Clearly, if $M_1 = M_2 = M_3$, then $G_c = G$. A core $G_c$ is {\em proper} if $G_c \not = G$. 

Proposition \ref{einstieg} follows from the strict version of Petersen's Theorem \cite{Petersen_1891} that
for every edge $e$ of a bridgeless cubic graph there is a 1-factor containing $e$. 

\begin{proposition} \label{einstieg}
Let $G$ be a cubic graph. If $G$ has a 1-factor, then it has a core. 
Furthermore, if $G$ is bridgeless, then for every $v \in V(G)$ there is a core $G_c$ such that $ v \not \in V(G_c)$.
\end{proposition}

\begin{lemma} \label{subgraph}
Let $G_c$ be a core of a cubic graph $G$, and $K_c$ be 
a component of $G_c$. Then ${\cal M}$ is a 1-factor of $G_c$, and \\
1) $K_c$ is either an even circuit or it is a subdivision of a cubic graph $K$, and \\
2) if $K_c$ is a subdivision of a cubic graph $K$, then $E(K_c) \cap E_3$ is a 1-factor of $K$,
and every edge of $K$ is subdivided by an even number of vertices.
\end{lemma} 
{\bf Proof.} Let $G_c$ be a core of $G$ with respect to three 1-factors $M_1$, $M_2$ and $M_3$. 
The set ${\cal M}$ is a matching 
in $G$. Hence, every vertex $v$ of $G_c$ is incident to at most one edge of ${\cal M}$.
Since $M_1$, $M_2$ and $M_3$ are 1-factors, $v$ cannot be incident to three edges of ${\cal U}$. Hence,
${\cal M}$ is a 1-factor of $G_c$. If $v$ is incident to an edge of $E_3$, then $d_{G_c}(v) = 3$, and if it is 
incident to an edge of $E_2$, then $d_{G_c}(v) = 2$. \\
1) Let $K_c$ be a component of $G_c$. If it has no trivalent vertices, then
$E(K_c) \cap {\cal M}$ is a 1-factor of $K_c$ and hence, $K_c$ is an even circuit, whose edges
are in ${\cal M}$ and ${\cal U}$, alternately. If $K_c$ contains trivalent vertices, then it is a subdivision of a 
cubic graph $K$. \\
2) The set $E(K_c) \cap E_3$ is a matching in $K_c$
that covers all trivalent vertices of $K_c$. Since $K$ is obtained from $K_c$ by suppressing the bivalent vertices and 
the edges of $E(K_c) \cap E_3$  are unchanged, it follows that $E(K_c) \cap E_3$ is a 1-factor $F$ of $K$.

Furthermore, every edge $e \in E(K) - F$ corresponds to a path in $K_c$ that starts and ends with an edge of ${\cal U}$. 
Hence, it is subdivided by an even number of vertices. \hfill $\square$

\begin{lemma} \label{cardinality_M}
Let $k \geq 0$. If $G_c$ is a $k$-core of a cubic graph $G$, 
then $|{\cal M}| = k - | E_3|$.
\end{lemma}
{\bf Proof.} 
Every vertex of an edge of $E_2$ has degree 1 in $G[{\cal U}]$, and 
every vertex of an edge of $ E_3$ has degree 2 in $G[{\cal U}]$. Hence,
$k = \frac{1}{2}[2 |E_2| + 4 | E_3| ] = |E_2| + 2|E_3| = |{\cal M}| + | E_3|$.
\hfill $\square$

The {\em dumbbell graph} is the unique cubic graph which is obtained from $K_2$ by adding a loop to each vertex.

\begin{lemma} \label{basics_gamma_3}
Let $k \geq 0$. If $G_c$ is a $k$-core of a connected cubic graph $G$, then \\
1) $k < 3$ if and only if $G$ is either the dumbbell graph or $G$ is 3-edge-colorable. \\
2) $|V(G_c)| = 2k - 2| E_3|$, and $|E(G_c)| = 2k - |E_3|$. \\
3) $girth(G_c) \leq 2k$. \\
4) $G_c$ has at most $2k/girth(G_c)$ components.
\end{lemma}
{\bf Proof.} Let $G_c$ be a $k$-core of $G$ with respect to three 1-factors $M_1$, $M_2$, $M_3$. 

1) Let $k <3$. If $G_c$ has a bridge $e$, then all edges which are adjacent to $e$ are elements of ${\cal U}$. Hence,
$G_c$ is the dumbbell graph, and consequently $G=G_c$.
If $G_c$ has no bridge, then there are $i$, $j$ such that $1 \leq i < j \leq 3$ and $M_i \cap M_j = \emptyset$. Hence,
$M_i \cup M_j$ is an even 2-factor of $G$ and therefore, $G$ is 3-edge-colorable. The other direction is trivial. \\
2) ${\cal M}$ is a 1-factor of $G_c$. Hence, $|V(G_c)| = 2k - 2| E_3|$ by Lemma \ref{cardinality_M}. Since
$|{\cal U}| = k$ and ${\cal M} \cap {\cal U} = \emptyset$ it follows that $|E(G_c)| = |{\cal M}| + |{\cal U}| = 2k - |E_3|$.\\
3, 4) Since $G_c$ has minimum degree 2 and at most $2k$ vertices, it follows that every component 
contains a circuit of length at most $2k$. Since ${\cal M}$ is a 1-factor of $G_c$, it follows that
each circuit contains at least $\frac{1}{2}girth(G_c)$ edges of ${\cal U}$. Hence, there are at most $\frac{2k}{girth(G_c)}$ pairwise disjoint 
circuits in $G_c$. 
\hfill $\square$

\begin{corollary} \label{girth_seymour}
Let $G$ be a loopless cubic graph. If $\mu_3(G) \not = 0$, then $\mu_3(G) \geq 3$ and $girth(G) \leq 2\mu_3(G)$.
\end{corollary}

\subsection{The conjecture of Fan and Raspaud}

If a core $G_c$ of a cubic graph is a cycle, then we say that $G_c$ is a {\em cyclic} core.
A cubic graph $G$ has three 1-factors such that $E_3 = \emptyset$ if and only if 
$G$ has a cyclic core.  Hence, 
Conjecture \ref{Fan_Raspaud} can be formulated as a conjecture on cores in bridgeless cubic graphs. 

\begin{conjecture} [Conj.~\ref{Fan_Raspaud}] \label{Fan_Raspaud_2}
Every bridgeless cubic graph has a cyclic core.
\end{conjecture}
Let $K_2^3$ be the unique cubic graph on two vertices which are connected by three edges. 

\begin{theorem} \label{girth_cyclic_core_thm} 
Let $k > 0$ and $G_c$ be a proper $k$-core of a cubic graph $G$ with $girth(G_c) \geq k$. Then \\
1) $G_c$ is a circuit of length $2k$, or \\
2) $k$ is even, $girth(G_c) = k$ and 
 $G_c$  is the disjoint union of two circuits of length $k$, or 
 $G_c$ is the graph $K_2^3$ with two of its edges subdivided by $k-2$ vertices.\\
In particular, if additionally $girth(G_c) > k$ or $k$ is odd, then $G_c$ is a circuit.
\end{theorem}

{\bf Proof.} Let $G_c$ be proper $k$-core of $G$ and $girth(G_c) \geq k$.
By Lemma \ref{basics_gamma_3}, $|E(G_c)| = 2k - | E_3|$ and hence, $|E(G_c)| \leq 2girth(G_c) -  | E_3|$ (*). 
Furthermore, $G_c$ has at most two components. 

If it has two components, then each of them contains a circuit. Hence, $|E(G_c)| \geq 2 girth(G_c) \geq 2k$. Thus, 
$| E_3|=0$ and $|E(G_c)|=2k$ and $G_c$ is the disjoint union of two circuits $C_1', C_2'$ of length $k$. 
By Lemma \ref{subgraph}, ${\cal M} \cap E(C_i')$ is a 1-factor of $C_i'$. Hence, $k$ is even.

Now suppose that $G_c$ is connected. If $|E_3|=0$, then it is a circuit of length $2k$. If $|E_3| > 0$,
then it contains at least two circuits and, by (*), any two circuits of $G_c$ intersect. 
Thus, $G_c$ is bridgeless. Let $e \in E_3$. 
Since $G_c -  E_3$ is a circuit $C_c$, it follows that $e$ is a chord of $C_c$. Thus,
there are two circuits $C_1$ and $C_2$ with $E(C_1) \cap E(C_2) = \{e\}$.
Hence,
$2girth(G_c) -  | E_3| \geq |E(G_c)| \geq |E(C_1) \cup E(C_2)| = |E(C_1)| + |E(C_2)| - 1 \geq 2 girth(G_c) -1$, and therefore,
$| E_3|=1$ and $|E(G_c)| = 2k-1$.
Thus, $G_c$ is a subdivision of $K_2^3$, where two edges are subdivided by $k-2$ vertices. Furthermore, $k$ is even by Lemma \ref{subgraph}. 
\hfill $\square$ 

Since $girth(G_c) \geq girth(G)$, Theorem \ref{girth_cyclic_core_thm} implies the following corollary.

\begin{corollary} \label{girth_cyclic_core_thm_corollary} 
Let $G$ be a cubic graph with $girth(G) \geq \mu_3(G)$. Then every $\mu_3(G)$-core is bipartite. In particular, if 
additionally $girth(G) > \mu_3(G)$
or $\mu_3(G)$ is odd, then every $\mu_3(G)$-core is a circuit. 
\end{corollary}

Let $G$ be a bridgeless cubic graph. The minimum number of odd circuits in a 2-factor of $G$ is the {\em oddness} of $G$, which is denoted by
$\omega(G)$. M\'a\v{c}ajov\'a and \v{S}koviera proved that 
Conjecture \ref{Fan_Raspaud} is true for bridgeless cubic graphs with oddness at most 2. 

\begin{theorem} [\cite{M_Skoviera_2009}] \label{cyclic_core_oddness_2}
Let $G$ be a bridgeless cubic graph. If $\omega(G) \leq 2$, then $G$ has a cyclic core.
\end{theorem}

We will use the following proposition for the proof of the next theorem. 

\begin{proposition} [\cite{Steffen_98}] \label{oddness_2_resistance}
Let $G$ be a bridgeless non-3-edge-colorable cubic graph. There is a proper 4-edge-coloring of $G$ with a color class that contains precisely 
two edges  if and only if $\omega(G)=2$.
\end{proposition}

\begin{theorem} \label{<6_cycle-core}
Let $G$ be a  simple bridgeless cubic graph. If $\mu_3(G) \leq 6$, then $G$ has a cyclic core. In particular, if $G$ is triangle-free
and $\mu_3(G) \leq 5$, then every $\mu_3(G)$-core is cyclic. 
\end{theorem}
{\bf Proof.} Let $G_c$ be a core of $G$ with respect to three 1-factors $M_1$, $M_2$, $M_3$, and $|E_0| = \mu_3(G)$. 
If $\mu_3(G) = 0$, then there is nothing to prove.
Since $G$ is bridgeless, it follows by Lemma \ref{basics_gamma_3} that $\mu_3(G) \geq 3$. Furthermore, $G$ has no loop and hence, 
if $\mu_3(G) = 3$, then $G_c$ is cyclic. 
We will use the following claim. 

\begin{claim} \label{circuit_alt}
If $C$ is a circuit whose edges are in ${\cal M}$ and ${\cal U}$ alternately, then 
for every $i \in \{1,2,3\}$ there is an edge $e \in {\cal M} \cap E(C)$ such that $e \not \in M_i$. 
\end{claim}

{\bf Proof.}
Suppose to the contrary that there is $i \in \{1,2,3\}$ such that $|{\cal M} \cap E(C)| = |M_i \cap E(C)|$. Let $i=3$. Then
$M_3' = (M_3 - E(C)) \cup ({\cal U} \cap E(C))$ is a 1-factor of $G$. With  
$M_1=M_1'$, $M_2=M_2'$ we get 
$|E_0'| = |E(G) - \bigcup_{i=1}^3 M_i'| < |E(G) - \bigcup_{i=1}^3 M_i| = \mu_3(G)$, a contradiction. $\square$

Claim \ref{circuit_alt} implies, that $G_c$ does not contain a circuit of length 4 whose edges are in ${\cal M}$ and ${\cal U}$ alternately.
Hence, every component of $G_c$ contains at least three edges of ${\cal U}$. If it has two components, then $\mu_3(G) = 6$ and 
each component is a circuit of length 6. Hence, $G_c$ is cyclic in this case.

Thus, we now assume that $G_c$ is connected.
Suppose to the contrary that there is an edge $e \in E_3$.

If no edge of $E_3$ is a bridge, then $G_c -E_3$ is a circuit. 
If $| E_3 | \in \{2,3\}$, then it follows that $G_c$ contains a circuit with edges in ${\cal M}$ and ${\cal U}$
alternately and there is an $i \in \{1,2,3\}$ such that $|{\cal M} \cap E(C)| = M_i \cap E(C)$. Hence, we obtain a contradiction with 
Claim \ref{circuit_alt}. 
Hence,  $| E_3|=1$ and $G_c$ is a subdivision of $K_2^3$ where two edges are subdivided by four vertices, that is $\mu_3(G) = 6$. 

(***) Consider 
$(\bigcup_{i=1}^3 M_i) - E(G_c)$ as proper 3-edge-coloring $\phi$ of $G-E(G_c)$. In any case $\phi$
can be extended to a proper 4-edge-coloring of $G$ which has a color class that contains precisely two edges. Now the result
follows with Proposition \ref{oddness_2_resistance} and Theorem \ref{cyclic_core_oddness_2}. 

For the remainder of the proof we suppose that $ E_3$ contains a bridge of $G_c$. If it has more than one bridge, then 
$G_c - E_3$ has at least three components. It is easy to see that $|{\cal U}| > 6$ in this case, contradicting the fact
that $\mu_3(G) \leq 6$. 

We now assume that $e$ is the only bridge of $G_c$. Then $e$ connects two disjoint circuits $C_1$ and $C_2$ which form a 2-factor of $G_c$.
Since $\mu_3(G) \leq 6$, it follows  that $|E(C_1)| + |E(C_2)| \in \{6,8,10\}$. 
If $\mu_3(G) \leq 5$, then one of these two circuits is a triangle. Hence, if $G$ is triangle-free, then we obtain a contradiction and the
statement for triangle-free graphs is proved. In the other cases we argue as in (***) that there is an appropriate proper
4-edge-coloring of $G$ such that $G$ has a cyclic core by  Proposition \ref{oddness_2_resistance} and Theorem \ref{cyclic_core_oddness_2}
\hfill $\square$

\subsection{Berge-Fulkerson Conjecture}

Following \cite{Seymour_79} we define a {\em $p$-tuple edge multicoloring} ($p > 1$) of a bridgeless cubic graph $G$ as a list of 
$3p$ 1-factors such that $E_p = E(G)$.

\begin{theorem} [\cite{Seymour_79}] \label{Paul_79}
Let $G$ be a bridgeless cubic graph which has no non-trivial 3-edge-cut. If $M$ is a 1-factor of $G$, then there are an integer $p > 1$ and a 
$p$-tuple edge multicoloring of $G$ using $M$. 
\end{theorem}

\begin{lemma} \label{seymour_lemma}
Let $G$ be a bridgeless cubic graph which has no non-trivial 3-edge-cut, 
$M$ a 1-factor of $G$ and $P$ a path of length 3. If $M$ contains no edge of $P$, then there is
a 1-factor $M'$ of $G$ that contains the two endedges of $P$.
\end{lemma}
{\bf Proof.} Let $P$ be a path with vertex set $\{v_1, \dots, v_4\}$ and edge set $\{ v_iv_{i+1} \colon 1 \leq  i \leq 3\}$, and let
$e= v_1v_2$ and $e'= v_3v_4$. We will show that there is a 1-factor $M'$ that contains $e$ and $e'$. 

Let
$f_2, f_3$ be the edges which are adjacent to $v_2$, $v_3$, respectively, and which are no edges of $P$. 
If  $f_2 = f_3$, (that is, $v_2$ and $v_3$ are connected by two edges) then every 1-factor that contains $e$ has to contain $e'$. 

If $f_2 \not= f_3$, then $f_2, f_3 \in M$. 
Theorem \ref{Paul_79} implies that there exist an integer $p > 1$ and a $p$-tuple edge multicoloring $\phi$ of $G$ using $M$.  Let
$M_1, \dots ,M_p$ be the $p$ 1-factors of $\phi$ that contain $e$.
If $M_i$ does not contain $e'$, then it contains $f_3$. Since 
$f_3 \in M$ and $M_i \not = M$, for all $i \in \{1, \dots,p\}$, there is 
$j \in \{1, \dots,p\}$ such that $f_3 \not \in M_j$. Hence, $M_j$ contains $e$ and $e'$. \hfill $\square$

\begin{theorem} \label{gamma_3_Fulkerson}
Let $G$ be a bridgeless cubic graph which has no non-trivial 3-edge-cut. If $\mu_3(G) \leq 4$, then $G$ has a Berge-cover.
\end{theorem}
{\bf Proof.} If $\mu_3(G) = 0$, then $G$ is 3-edge-colorable, and it has a Berge-cover. 
Let $M_1$, $M_2$, $M_3$ be three 1-factors of $G$ such that $\mu_3(G) = |{\cal U}|$, and $G_c$ be the induced core.  
Then $\mu_3(G) \geq 3$, by Corollary 
\ref{girth_seymour}. Since $G$ is triangle-free, it follows with Theorem \ref{<6_cycle-core} that $G_c$ is cyclic.
Because of the minimality of $\mu_3(G)$ it follows that $G_c$ is connected. Hence, 
the edges of ${\cal U}$ can be paired into at most two pairs
(one pair and a single edge if $\mu_3(G) = 3$), such that the edges of a pair are connected by an edge of ${\cal M}$. 
Lemma \ref{seymour_lemma} implies that there are two 1-factors $M_4$ and $M_5$ such that 
$\bigcup_{i=1}^5 M_i = E(G)$.
\hfill $\square$

\subsection{Short cycle covers}

\begin{theorem} \label{scc_core_allgemein}
Let $k,l,t$ be non-negative integers, and $G$ be a cubic graph. If $G$ has a $k$-core which has a $l$-cycle cover ${\cal C}_c$ of length $t$, 
then $G$ has a $(l+2)$-cycle cover ${\cal C}$ of length at most $\frac{4}{3}(|E(G)|-k) + t$. Furthermore, ${\cal C}$  is even if and only if 
${\cal C}_c$  is even. 
\end{theorem}
{\bf Proof.} 
Let $G_c$ be a $k$-core of $G$ with respect to three 1-factors $M_1$, $M_2$, $M_3$. 
Then $\{ M_1 \oplus M_2, M_1 \oplus M_3, M_2 \oplus M_3\}$ are three even cycles which together cover, and only cover, each 
edge of $E_1 \cup E_2$ precisely twice. Thus, there are two of these three cycles covering $E_1 \cup E_2$ with total length at most
$\frac{4}{3}|E_1 \cup E_2| \leq \frac{4}{3}(|E(G)| - k)$. These two cycles together with a cycle cover of $G_c$ gives a cycle cover
of $G$ with length at most $ \frac{4}{3}(|E(G)| - k) + t$. The cycle cover is even if and only if the cycle cover of $G_c$ is even.
\hfill $\square$

The following result of Alon and Tarsi [1] and Bermond, Jackson and Jaeger [2]  is the best known general result on the length of cycle covers.

\begin{theorem} [\cite{Alon_Tarsi, Bermond_5/3}] \label{5/3}
Every bridgeless graph $G$ has a 3-cycle cover of length at most $\frac{5}{3}|E(G)|$.
\end{theorem}

\begin{theorem} \label{scc_bridgeless_core}
Let $k \geq 0$, and $G$ be a cubic graph. If $G$ has a bridgeless $k$-core,
then $G$ has a cycle cover of length at most $\frac{4}{3}|E(G)| + 2k$.
\end{theorem}
{\bf Proof.} Let $G_c$ be a bridgeless $k$-core of $G$. By Theorem \ref{5/3}, $G_c$ has 
a cycle cover of length at most $\frac{5}{3}|E(G_c)|$. By Lemma \ref{basics_gamma_3} we have 
$|E(G_c)| = 2k - | E_3|$, and hence, it follows with Theorem \ref{scc_core_allgemein} that $G$ has a 
cycle cover of length at most $\frac{4}{3}|E(G)| + 2k$.\hfill $\square$

We are going to prove better bounds for the length of cycle covers of a cubic graphs which have a bipartite core. 
First we show that bipartite cores are bridgeless.

\begin{theorem} \label{bipartite_bridgeless}
Let $G_c$ be a core of a cubic graph. If $G_c$ is bipartite, then $G_c$ is bridgeless.
\end{theorem}
{\bf Proof.}
If $G_c$ is not bridgeless, then it has a component $K_c$ that contains a bridge.
Furthermore, $K_c$ has a bridge $e$ such that one component of $K_c-e$ is 2-edge-connected. Let 
$K_c'$ be a 2-edge-connected component of $K_c - e$.
Since $e \in E_3$ and, by Lemma \ref{subgraph}, ${\cal M} \cap E(K_c)$ is a 1-factor of $K_c$,
it follows that $|V(K_c')|$ is odd. Furthermore, if we remove all edges of 
$ E_3$ from $K_c'$, then we obtain a set of circuits. 
Hence, $G_c$ contains an odd circuit. Therefore, it is not bipartite. 
\hfill $\square$

\begin{theorem} \label{scc_strong_core}
Let $k \geq 0$, and $G$ be a cubic graph. If $G$ has a bipartite $k$-core $G_c$,
then $G$ has an even 4-cycle cover of length at most $\frac{4}{3}|E(G)| + \frac{2}{3}k$.
In particular, if  $G_c$ is cyclic, then $G$ has an even 3-cycle cover of length at most $\frac{4}{3}|E(G)| + \frac{2}{3}k$.
\end{theorem}
{\bf Proof.} 
Let $G_c$ be a bipartite $k$-core of $G$. 
Let $K_c$ be a component of $G_c$, and $k' = |E(K_c) \cap {\cal U}|$. 
It suffices to prove that every component  $K_c$ has an even cycle cover of length $2k'$.
Then it follows that $G_c$ has an even cycle cover of length $2k$. Hence,
$G$ has an even cycle cover of length at most $\frac{4}{3}|E(G)| + \frac{2}{3}k$ by Theorem \ref{scc_core_allgemein}.

Clearly, if $K_c$ is a circuit, then it has a 1-cycle cover of length $2k'$. Thus, the statement is true for $G_c$ is cyclic.

If $K_c$ is a component of $G_c$ which is not a circuit, then $E(K_c) \cap E_3 \not = \emptyset$, and by Theorem \ref{bipartite_bridgeless},
$G_c$ is 2-edge-connected. Furthermore, $K_c$ is a subdivision 
of a cubic graph $H_c$. Let $E^* = E(K_c) \cap E_3$.
Since the two vertices which are incident to an edge of $E^*$ have degree 3 in $K_c$,
it follows that $K_c - E^*$ is a cycle. Since $K_c$ is bipartite, it is an even cycle which has a proper 2-edge-coloring. Thus, 
$E^*$ is a color class of a proper 3-edge-coloring $\phi$ of $K_c$.
By Lemma \ref{subgraph}, each edge of $H_c$ is subdivided by an even number of vertices. 
Hence, $\phi$ induces a proper 3-edge-coloring $\phi'$ on $H_c$, where $E^*$ is a color class. 
Let ${\cal C}_{H_c}$ be a cycle cover of $H_c$ of length $\frac{4}{3}|E(H_c)|$, which is induced by $\phi'$ and which uses 
the edges of $E^*$ twice.

Now, ${\cal C}_{H_c}$  induces a 2-cycle cover ${\cal C}_{K_c}$ of $K_c$.
The length of ${\cal C}_{K_c}$ is $|E(K_c)| + |E^*|$. 
By Lemma \ref{basics_gamma_3} we have
$|E(K_c)| = 2k' - |E^*|$. Hence, the length of ${\cal C}_{K_c}$ is $2k'$.

It remains to show that  ${\cal C}_{K_c}$ is an even cycle.
Every edge of $E^*$ is contained in precisely 
two cycles of ${\cal C}_{K_c}$ and all the other in precisely one. Let $v$ be a vertex which is incident to two edges of ${\cal U}$. Then
$d_{G_c}(v)=3$, and $v$ is incident to an edge of $E^*$. Hence, every circuit $C$ of ${\cal C}_{K_c}$ does not contain any two
consecutive edges of ${\cal U}$. Since ${\cal M}$ is a 1-factor of $G_c$, it follows that the edges of $C$ are in ${\cal M}$ and ${\cal U}$
alternately. Hence, $C$ has even length, and ${\cal C}_{K_c}$ is an even 2-cycle cover of  $K_c$.
\hfill $\square$

\begin{corollary} 
Let $G$ be a triangle-free bridgeless cubic graph.
 If $\mu_3(G) \leq 5$, then 
$G$ has an even 3-cycle cover of length at most $\frac{4}{3}|E(G)| + 2$.
\end{corollary} 
{\bf Proof.} If $\mu_3(G) \leq 5$, then Theorem \ref{<6_cycle-core} implies that the induced core is cyclic. 
If $\mu_3(G) = 5$, then $\frac{4}{3}|E(G)| + \lfloor \frac{2}{3} \mu_3(G) \rfloor $ is odd. Hence, the 
result follows with Theorem \ref{scc_strong_core}. \hfill $\square$

In \cite{Fan_Raspaud_94} it is proved that if a cubic graph $G$ has a cyclic core, then 
 it has a 3-cycle cover of length at most $\frac{14}{9}|E(G)|$. This result is improved to smaller than $\frac{14}{9}|E(G)|$ in \cite{Macajova_etal_2011}. We additionally deduce that there is an even 3-cycle cover of length smaller than  $\frac{14}{9}|E(G)|$.

\begin{corollary} Let $G$ be a cubic graph. If $G$ has a cyclic core, then $G$ has an even 3-cycle cover of length smaller than
$\frac{14}{9}|E(G)|$.
\end{corollary}
 {\bf Proof.} If $G$ is 3-edge-colorable, then the statement is true. Let $G_c$ be a cyclic $k$-core of $G$. Then $G_c$ is 2-regular and it has 
$2k$ edges. Hence, $2k \leq \frac{2}{3}|E(G)|$ and therefore, $k \leq \frac{1}{3}|E(G)|$. If $k = \frac{1}{3}|E(G)|$, then $G_c$ is an even
2-factor of $G$ and hence, $G$ is 3-edge-colorable. Thus, $k < \frac{1}{3}|E(G)|$ and Theorem 
\ref{scc_strong_core} implies that $G$ has an even 3-cycle cover of length smaller than $\frac{14}{9}|E(G)|$. \hfill $\square$

In \cite{Fan_Raspaud_94} it is proved that if a cubic graph $G$ has a Fulkerson-cover, then 
it has a 3-cycle cover of length at most $\frac{22}{15}|E(G)|$. This bound is best possible for 3-cycle covers of bridgeless cubic graphs, 
since it is attained by the Petersen graph \cite{Macajova_etal_2011}. We additionally show that there exists such a cycle cover which is even.

\begin{corollary} \label{Fulkerson_scc}
Let $G$ be a cubic graph which has a Fulkerson-cover. \\
1) Then $G$ has an even 3-cycle cover of length at most $\frac{22}{15}|E(G)|$.\\
2) If $|V(G)| \not \equiv 0 \bmod 10 $, then $G$ has an even 3-cycle cover of length smaller than  $\frac{22}{15}|E(G)|$.
\end{corollary}
{\bf Proof.} 1) Let $M_1, \dots, M_6$ be the six 1-factors of a Fulkerson-cover of $G$.
Since ${ 6 \choose 2 } = 15$, there are two 1 factors, say $M_1, M_2$, such that $|M_1 \cap M_2| \leq \frac{1}{15}|E(G)|$.
We claim that there is $i \in \{3, \dots ,6\}$ such that $|M_1 \cap M_2| + |M_1 \cap M_i| + |M_2 \cap M_i| \leq \frac{1}{5}|E(G)|$. 
Suppose to the contrary that this is not true. Then 
$ \sum_{i=3}^6  (|M_1 \cap M_2| + |M_1 \cap M_i| + |M_2 \cap M_i|) > \frac{4}{5}|E(G)|$. We have 
$ \sum_{i=3}^6  (|M_1 \cap M_2| + |M_1 \cap M_i| + |M_2 \cap M_i|) = \frac{2}{3}|E(G)| + 2|M_1 \cap M_2|$ and hence,
$|M_1 \cap M_2| > \frac{1}{15}|E(G)|$, which contradicts our choice of $M_1$ and $M_2$. 

Let $i = 3$ and $|M_1 \cap M_2| + |M_1 \cap M_3| + |M_2 \cap M_3| \leq \frac{1}{5}|E(G)|$. Since 
every edge is contained in precisely two 1-factors, the $k$-core with respect to $M_4$, $M_5$ and $M_6$ is
cyclic and $k \leq \frac{1}{5}|E(G)|$. Theorem \ref{scc_strong_core} implies that $G$ has an even 3-cycle cover of length
at most  $\frac{22}{15}|E(G)|$. 

2) If $|V(G)| \not \equiv 0 \bmod 10$, then $|E(G)| \not \equiv 0 \bmod 15$, and 
we deduce as above that $G$ has cyclic $k$-core and 
$k < \frac{1}{5}|E(G)|$. Then the statement follows with Theorem \ref{scc_strong_core}
\hfill $\square$

Let $G$ be a cubic graph which has a 1-factor and consequently a core. If $G$ has a bridge, then every core of $G$ has a bridge. 
We conjecture that the opposite direction of that statement is true as well and  propose two conjectures.

\begin{conjecture} \label{Conj_bridgeless}
Every bridgeless cubic graph has a proper bridgeless core. 
\end{conjecture}

\begin{conjecture} \label{Conj_bipartite}
Every bridgeless cubic graph has a proper bipartite core. 
\end{conjecture}

Conjecture \ref{Fan_Raspaud} implies Conjecture \ref{Conj_bipartite}, which implies Conjecture \ref{Conj_bridgeless},
by Theorem \ref{bipartite_bridgeless}.

\section[]{Four 1-factors}  \label{gamma_4}

Let ${\cal C}$ be a cycle cover of a graph $G$. For $e \in E(G)$, let $ced_{{\cal C}}(e) = |\{C \colon C \in {\cal C} \mbox{ and }e \in E(C)\}|$,
and $\max \{ced_{{\cal C}}(e) \colon e \in E(G) \}$ be the {\em edge-depth of} ${\cal C}$, which is denoted by $ced_{{\cal C}}(G)$.
Zhang 
conjectured that every bridgeless graph has a shortest cycle cover of at most four cycles (Conjecture 8.11.5 in \cite{Zhang_book}), and 
that every 3-edge-connected graph has a shortest cycle cover ${\cal C}$ such that $ced_{{\cal C}}(G) \leq 2$ (Conjecture 8.11.6 in \cite{Zhang_book}).
The next theorem shows that we get the optimal bound for the length of a cycle cover if $\mu_4(G)=0$, and that these graphs have 
a 5-cycle double cover. It also shows that these graphs satisfy the aforementioned conjectures of Zhang as well.

\begin{theorem} \label{gamma_4_result}
Let $G$ be a cubic graph. If $\mu_4(G) = 0$, then \\
1) $G$ has an even 4-cycle cover ${\cal C}$ of length $\frac{4}{3}|E(G)|$, and $ced_{{\cal C}}(G) \leq 2$.\\
2) $G$ has a 5-cycle double cover.
\end{theorem}
{\bf Proof.} The statements are true for 3-edge-colorable cubic graphs. We assume that $G$ is not 3-edge-colorable.
Let $M_1, \dots, M_4$  be four 1-factors of $G$ with $\mu_4(G)=0$. Since $G$ is cubic it follows that $E_3 = E_4 = \emptyset$,
$E_2$ is a 1-factor and $E_1$ is a 2-factor of $G$. For $i \in \{1, \dots, 4\}$ let $F_i = E_2 \oplus M_i$. 
Then, ${\cal F} = \{F_1, F_2, F_3, F_4\}$
is an even 4-cycle cover that covers each edge of $E_1$ once and each edge of $E_2$ twice. Hence, it has length 
$|E_1| + 2 |E_2| = \frac{4}{3}|E(G)|$. Since $E_1$ is a 2-factor it follows that ${\cal F} \cup E_1$ is a 5-cycle double cover of $G$.
\hfill $\square$

Theorem \ref{gamma_4_result} 2) was proved by Hou, Lai and Zhang \cite{Hou_etal_2012} independently. 

Following Zhang \cite{Zhang_book} we say that the Chinese postman problem is equivalent to the shortest cycle cover problem,
if the shortest length of a closed trail that covers all edges of $G$ is equal to the length of a shortest cycle cover. This is certainly true
for cubic graphs $G$ that have a cycle cover of length $\frac{4}{3} |E(G)|$. Zhang  asked the following question (Problem 8.11.4 in \cite{Zhang_book}):
Let $h \geq 5$ and $G$ be a 3-edge-connected, cyclically $h$-edge-connected graph. If the Chinese Postman problem and the 
shortest cycle cover problem are equivalent for $G$, does $G$ admit a nowhere-zero 4-flow?  
The answer to this question is negative since for $h \in \{5,6\}$ there are cyclically $h$-edge connected snarks with $\mu_4(G) = 0$. 
It is known that $\mu_4(G) = 0$ if $G$ is a flower snark or a Goldberg snark, see \cite{Fouquet_09}.

We now study the case, when the union of four 1-factors does not cover all edges of a cubic graph $G$. 

\begin{lemma} \label{4cyclecover}
Let $G$ be a cubic graph that has four 1-factors $M_1, \dots, M_4$ with $|E(G) - \bigcup_{i=1}^4 M_i| = k \geq 0$.
If $E_4= \emptyset$, then $G$ has a 4-cycle cover of length $\frac{4}{3}|E(G)| + 4k$.
\end{lemma}
{\bf Proof.} If $G$ is 3-edge-colorable or $k=0$, then the statement is true by Theorem \ref{gamma_4_result}.
Let $G$ be not 3-edge-colorable, and $k > 0$. Since $E_4 = \emptyset$, it follows that $G$ is bridgeless. 

For $ i \in \{1, \dots ,4\}$ and $ j \in \{1,2,3\}$ let $M_i^j = M_i \cap E_j$ and $\bar{M_i}^j = (E(G) - M_i) \cap E_j$. 
For $i \in \{1, \dots,4\}$, let ${\cal C}_i' = M_i^1 \cup \bar{M_i}^2 \cup M_i^3 \cup {\cal U}$. 
Every vertex of an edge of $M_i^1$ is incident either 
to an edge of $\bar{M_i}^2$ and to an edge of $\bar{M_i}^1$, or
to an edge of $\bar{M_i}^3$ and to an edge of ${\cal U}$. 
Every vertex of an edge of $\bar{M_i}^2$ is incident either 
to an edge of $M_i^1$ and to an edge of $\bar{M_i}^1$, or
to an edge of $M_i^2$ and to an edge of ${\cal U}$.  
Every vertex of an edge of $M_i^3$ is incident to an edge of $\bar{M_i}^1$ and an edge of ${\cal U}$.  
Every vertex of an edge of  ${\cal U}$  is incident either 
to an edge of $M_i^1$ and to an edge of $\bar{M_i}^3$, or
to an edge of $M_i^2$ and to an edge of $\bar{M_i}^2$, or 
to an edge of $M_i^3$ and to an edge of $\bar{M_i}^1$.
Hence, every vertex of $G[{\cal C}_i']$ is adjacent to precisely two edges of 
${\cal C}_i'$; that is, ${\cal C}_i'$ is a cycle. Note that the edges of $M_i^2$ (and hence, the vertices as well) are not
in $G[{\cal C}_i']$.

Thus, $\bigcup_{i=1}^{4} {\cal C}_i' = E(G)$. Let ${\cal C}' = \{ {\cal C}_1', {\cal C}_2', {\cal C}_3', {\cal C}_4' \}$. Then ${\cal C}'$ 
is a 4-cycle cover of $G$. Each $e \in E(G)$ is either an element of ${\cal U}$ or there are  
$i \in \{1, \dots,4\}$ and $j \in \{1,2,3\}$ such that $e \in M_i^j$. If $e \in M_i^j$, then it is contained in precisely $j$ 
cycles and if $e \in {\cal U}$, then it is contained in all four cycles. Hence, the length of ${\cal C}'$ 
is $ \frac{4}{3}|E(G)| + 4k$.
\hfill $\square$

\begin{theorem} \label{gamma_4_corollary}
Let $G$ be a loopless cubic graph. If $\mu_4(G) \in \{0,1,2,3\}$, then $G$ has a 4-cycle cover of length $\frac{4}{3}|E(G)| + 4 \mu_4(G)$.
\end{theorem}
{\bf Proof.} 
Since $\mu_4(G) \leq 3$, it follows that $E_4 = \emptyset$. The 
result follows with Lemma \ref{4cyclecover}.  \hfill $\square$ \\
The bound of Corollary \ref{gamma_4_corollary} is attained by the Petersen graph $P$ with  $\mu_4(P)=1$.

\section{Remark on hypohamiltonian snarks}

A graph $G$ is {\em hypohamiltonian} if it is not hamiltonian but $G-v$ is hamiltonian for every vertex $v \in V(G)$. Non 3-edge-colorable, cubic hypohamiltonian graphs are cyclically 4-edge-connected and have girth at least 5, and there are cyclically
6-edge-connected hypohamiltonian snarks with girth 6, see \cite{M_Skoviera_2007}. Since
hamiltonian cubic graphs are 3-edge-colorable, and $G-v$ is not 3-edge-colorable for every snark $G$, hypohamiltonian snarks could be considered as being
closest to being 3-edge-colorable. Hypohamiltonian snarks have a proper 4-edge-coloring with a color class of cardinality 2. Hence,
it follows with Proposition \ref{oddness_2_resistance} and Theorem \ref{cyclic_core_oddness_2} that they have a cyclic core.

Corollary \ref{girth_seymour} implies, that if $G$ is a cubic graph and $\mu_3(G)=3$, then 
$G$ has girth at most 6.
It is easy to see that $\mu_3(G) = 3$, if $G$ is the Petersen graph or a flower snark, which are hypohamiltonian snarks. 
Jaeger and Swart \cite{Jaeger_Swart_80} conjectured that (1) the girth and (2) the cyclic connectivity of a snark is at most 6. The first conjecture
is disproved by Kochol \cite{Kochol_96} and the second is still open. 
We believe that both statements of Jaeger and Swart are true for hypohamiltonian snarks.  

\begin{conjecture} \label{conj_hypo}
Let $G$ be a snark. If $G$ is hypohamiltonian, then $\mu_3(G) = 3$. 
\end{conjecture} 

H\"aggkvist \cite{Haeggkvist_2003} proposed to prove the Fulkerson conjecture for hypohamiltonian graphs, 
which might be easier to prove than the general case. 
By  Theorem \ref{gamma_3_Fulkerson}, Conjecture \ref{conj_hypo} implies 
that  hypohamiltonian snarks  have a Berge cover, and together with 
Theorem \ref{scc_strong_core} it implies that they have an even 3-cycle cover of length at most $\frac{4}{3}|E(G)| + 2$. 

\subsubsection*{Acknowledgement} I thank the anonymous referees for their helpful comments and suggestions. In particular, the short proof of Theorem \ref{gamma_4_result} was given by one of them.

\end{document}